\documentclass{amsart}

\usepackage{amssymb}

\theoremstyle{plain}

\numberwithin{equation}{section}

\newtheorem*{teorema}{Theorem A}
\newtheorem*{teoremab}{Theorem B}
\newtheorem*{teoremac}{Theorem C}

\newtheorem{lem}[equation]{Lemma}

\newcommand{\Irr}{\operatorname{Irr}}
\newcommand{\Lin}{\operatorname{Lin}}

\newcommand{\Ker}{\operatorname{Ker}}

\newcommand{\SL}{\operatorname{SL}}

\newcommand{\V}{\operatorname{V}}

\theoremstyle{definition}

\begin{document}	
\title{Products of characters and finite $p$-groups III}

\author{Edith Adan-Bante}

\address{University of Southern Mississippi Gulf Coast, 730 East Beach Boulevard,
 Long Beach MS 39560}

\email{Edith.Bante@usm.edu}


\subjclass{20c15}

\date{2003}

\begin{abstract} Let $G$ be a finite $p$-group and $\chi,\psi$ be irreducible 
characters of $G$. We study the character $\chi\psi$ when 
$\chi\psi$ has at most $p-1$ distinct irreducible constituents. 
\end{abstract}

\maketitle

\begin{section}{Introduction}

Let $G$ be a finite $p$-group. Denote by $\Irr(G)$ the set of 
irreducible complex characters of $G$. Let $\chi, \psi \in \Irr(G)$. Then the product
of $\chi\psi$ can be written as
$$
\chi\psi= \sum_{i=1}^n a_i\theta_i
$$
\noindent where $\theta_i\in \Irr(G)$ and $a_i=[\chi\psi,\theta]>0$. 
Set $\eta(\chi\psi)=n$. So $\eta(\chi\psi)$ is 
 the number of distinct irreducible constituents of the product $\chi\psi$.
The purpose of this
note is to study the case when the product $\chi\psi$ has ``few" irreducible
constituents, namely when $\eta(\chi\psi)<p$. 

If $\chi$ is a character of $G$, denote by $\V(\chi)=<g\in G\mid \chi(g)\neq 0>$.
So $\V(\chi)$ is  the smallest subgroup of $G$ such that $\chi$ vanishes on 
$G \setminus \V(\chi)$. Through this work,
we use the notation of \cite{isaacs}. In addition, we are going to 
denote by $\Lin(G)=\{\lambda\in \Irr(G) \mid \lambda(1)=1\}$ the set of
linear characters. 
The main results are 
	
\begin{teorema} 
Let $G$ be a finite $p$-group and $\chi, \psi \in \Irr(G)$.
Assume that $\eta(\chi\psi)<p$.
 Let $\theta\in \Irr(G)$ be a constituent
of $\chi\psi$.   Then 

(i) ${\bf Z}(\chi\psi)={\bf Z}(\theta)$.

(ii) $\V(\chi\theta)\geq \V(\theta)$. Therefore $\V(\chi)\cap\V(\psi)\geq \V(\theta)$.
\end{teorema}

\begin{teoremab} Let $G$ be a finite $p$-group and $\chi, \psi \in \Irr(G)$.
Assume that $\eta(\chi\psi)<p$.  Let $N$ be normal subgroup of $G$ and 
$\alpha, \gamma \in \Irr(H)$.
If $\alpha^G=e\chi$, for some integer $e$, and  $[\gamma, (\chi\psi)_N]\neq 0$, 
then $\gamma^G$ is a multiple
of an irreducible. In particular, if $|G:N|=p$ then $\gamma^G \in \Irr(G)$.
\end{teoremab}

\begin{teoremac}
Let $p$ be a finite $p$-group and $\chi\in \Irr(G)$. 

(i) If $\chi\neq 1_G$ then  $[\chi^2,\chi]=0$

(ii) If $p\neq 2$  and $\chi(1)>1$ then $[\chi^2,\lambda]=0$ for all $\lambda \in \Lin(G)$. 

(iii) Assume also that either $p\neq 2$ or $\eta(\chi^2)<p$. Then
there exists a subgroup $H$ of $G$ and $\alpha \in \Lin(H)$ such that 
$\alpha^G=\chi$ and $(\alpha^2)^G \in \Irr(G)$. Thus $\chi^2$ has an irreducible
constituent of degree $\chi(1)$. 
\end{teoremac}

Let $G$ be a $p$-group and  $\chi\in \Irr(G)$. Assume that $\chi(1)=p^n$ with $n\geq 1$.
Denote by  $\overline{\chi}$ the complex conjugate 
of $\chi$. In \cite{edith1} it is proved that $\eta(\chi\overline\chi)\geq 2n(p-1)+1$.
Thus $\eta(\chi\overline{\chi})>p$. 
 We can check that 
 the principal character $1_G\in \Irr(G)$ is a constituent of 
$\chi\overline{\chi}$. Thus Theorem A (i) and (ii) may not hold
true without the condition of $\eta(\chi\psi)<p$. 
Also, since $\chi(1)>1$ and $G$ is a $p$-group, there exist a 
subgroup $N$ of $G$ and $\alpha\in \Irr(N)$
such that $|G:N|=p$ and $\alpha^G=\chi$. Observe that $[1_N, (\chi\overline{\chi})_N]\neq 0$
and $\eta(1_N^G)=p$. Thus Theorem B may not hold true without the condition of
$\eta(\chi\psi)<p$. 

Consider the group $\SL(2,3)$ and the character $\chi\in \Irr(\SL(2,3))$ such 
that $\chi(1)=3$. 
We can check that
 $[\chi^2,\chi]=2$. Thus Theorem C (i) may not hold true without the condition 
 of $G$ being a $p$-group.   
Let $G$ be the dihedral group of order $8$ and $\chi\in \Irr(G)$ with $\chi(1)=2$.
We can check that 
$$\chi^2=\sum_{\lambda \in \Lin(G)} \lambda.$$  
Thus Theorem C (ii) and (iii) may not hold true without the additional hypotheses. 


\end{section}

\begin{section}{Proofs}
Let $H$ be a subgroup of $G$ and $\lambda\in \Irr(H)$. Denote by
$\Irr(G\mid \lambda)=\{\chi\in \Irr(G)\mid [\chi_H,\lambda]\neq 0\}$
the set of irreducible characters of $G$ lying above $\lambda$.

\begin{proof}[Proof of Theorem A]
{\bf (i)} Observe that ${\bf Z}(\chi\psi)\leq {\bf Z}(\theta)$. Set
$Z={\bf Z}(\chi\psi)$. So $(\chi\psi)_Z=\chi(1)\psi(1) \gamma$
for some $\gamma \in \Lin(Z)$.

Suppose that 
$Z< {\bf Z}(\theta)$. Let $Y/Z$ be chief factor of $G$ such that
$Y\leq {\bf Z}(\theta)$. Since $Y/Z$ is a chief factor of a $p$-group, it is cyclic.
Thus the character $\gamma\in \Lin(Z)$ extends to $Y$. We can check that  
$$(\chi\psi)_Y= \frac{\chi(1)\psi(1)}{p}\sum_{\delta \in \Lin(Y\mid\gamma)} \delta.$$ 
If there exists some $G$-invariant $\delta$ in $\Lin(Y\mid \gamma)$, then
all the elements in 
$\Lin(Y\mid \gamma)$ are $G$-invariant since $G$ is a $p$-group and $Y/Z$ is a chief 
factor. But then $\chi\psi$ must have at least
$p$ distinct irreducible constituents since $\theta_Y=\theta(1)\delta$ for 
some $\delta \in \Lin(Y\mid \gamma)$.
Therefore  the set $\Lin(Y\mid \gamma)$ forms a $G$-orbit. 
In particular ${\bf Z}(\theta)$ can not contain $Y$. We conclude that $Z={\bf Z}(\theta)$.

{\bf (ii)} 
Suppose that $\V(\chi\psi)\cap \V(\theta)< \V(\theta)$. 
Set $V=\V(\theta)$. We are going to conclude that necessarily we have that
$\chi\psi$ has at least $p$ distinct irreducible constituents.

Let $V/W$ be a chief factor
of $G$ such that $W\geq \V(\chi\psi)\cap V$. Observe that
if $g\in \V(\theta)\setminus W$, then $g\not\in \V(\chi\psi)$.
Therefore, by definition of $\V(\chi\psi)$, we have that
$(\chi\psi)(g)=0$ for all $g\in \V(\theta)\setminus W$. In particular,
for all $\gamma \in \Lin(V/ W)$ 
 we have that
 \begin{equation}\label{counting}
(\chi\psi)_V \gamma=(\chi\psi)_V.
\end{equation}	
	Let $v\in V\setminus W$ such that $\theta(v)\neq 0$. Observe that 
such an element exist because otherwise $W=\V(\theta)=V$.
Let $\sigma\in \Irr(V)$ be a constituent of $(\theta)_V$
and $\{\sigma^g\mid g \in T\}$ be the $G$-orbit of $\sigma$ in $G$.
By Clifford Theory we have that  $\theta_V= e \sum_{g\in T} \sigma^g$,
for some positive integer $e$. Thus
\begin{equation}\label{count}
 \sum_{g\in T}\sigma^g(v)=\frac{\theta(v)}{e}\neq 0. 
\end{equation} 
 Since
$G$ is a $p$-group and $V/W$ is a chief factor of $G$, we have that
$G$ acts trivially on $V/W$. Therefore the set of 
characters $\Lin(V/W)$ are $G$-invariant. Thus the $G$-orbit of $\sigma\gamma$
is the set $\{\sigma^g\gamma\mid g\in T\}$. By \eqref{counting}
$\sigma\gamma$ is an irreducible  constituent of $(\chi\psi)_V$.
Thus, there exists some character 
$\theta_{\gamma} \in \Irr(G)$ such that $[(\theta_{\gamma})_V, \sigma\gamma]\neq 0$.  
By Clifford theory we have that $(\theta_{\gamma})_V= f \sum_{g\in T} \sigma^g\gamma$
for some positive integer $f$.
Therefore 
\begin{equation}
\frac{\theta(v)}{e}\gamma(v)= [\sum_{g\in T} \sigma^g(v)] 
\gamma(v)=\sum_{g\in T} (\sigma^g\gamma)(v)=\frac{\theta_{\gamma}(v)}{f}.
\end{equation}
We conclude that if $\gamma(v)\neq 1$, then $\theta(v)\neq \theta_{\gamma}(v)$ and
therefore $\theta\neq \theta_{\gamma}$. Similarly we can check that
if $\delta\in \Lin(V/W)$ and $\delta\neq \gamma$, then there exists a constituent
$\theta_{\delta} \in \Irr(G)$ of the product $\chi\psi$ such that $[\delta, (\theta_{\delta})_V]\neq 0$ and $\theta_{\gamma}\neq \theta_{\delta}$.
Since $\Lin(V/W)$ has $p$ distinct irreducible  constituents, then  
 $\chi\psi$ has at least $p$ distinct irreducible constituents. 
A contradiction with our hypothesis. Therefore $\V(\theta)\leq \V(\chi\psi)$.
\end{proof}


\begin{lem}\label{countingabove}
 Let $G$ be a finite $p$-group and $N$ be a normal subgroup
of $G$. Let $\phi\in \Irr(N)$.
Then the set 
$\,\Irr(\, G\mid \phi\, )$ of all $\chi \in \Irr(G)$ 
lying over $\phi$ has either one 
or at     
least $p$ members.
\end{lem}
\begin{proof}
 Let $G_{\phi}$ be the stabilizer of $\phi$ in $G$. By Clifford theory we have that  $\eta(\phi^{G_{\gamma}})=\eta(\phi^{G})$. If $|G_{\phi}|<|G|$, by induction we have that
 either 
 $\eta(\phi^{G_{\gamma}})=1$ or $\eta(\phi^{G_{\gamma}})\geq p$ and the result holds.
We may assume that  $|G_{\phi}|=|G|$.
 In Lemma 4.1 \cite{edith}, it is  proved that the result holds if $\phi$ is a $G$-invariant character.
\end{proof}
%

\begin{proof}[Proof of Theorem B]
Since $N$ is normal in $G$ and $\gamma^G=\chi$, the irreducible constituents of
$\chi_N$ are of the form $\alpha^g$, for some $g\in G$, and $(\alpha^g)^G=e\chi$.
Since $[\gamma, (\chi\psi)_N]\neq 0$, there exist some $g\in G$ and some $\beta \in \Irr(N)$ 
such that $[\alpha^g\beta, \gamma]\neq 0$. 
By Exercise 5.3 of \cite{isaacs} we have that $(\alpha^g\psi_N)=(\alpha^g)^G\psi=e\chi\psi$.
Thus the irreducible constituents of $(\alpha^g\psi_N)^G$ are irreducible constituents of
$\chi\psi$. In particular, the irreducible constituents of $(\alpha^g\beta)^G$ are 
irreducible constituents of $\chi\psi$. Thus the irreducible constituents of $\gamma^G$
are irreducible constituents of $\chi\psi$. By Lemma \ref{countingabove}, either 
$\eta(\gamma^G)=1$ or $\eta(\gamma^G)\geq p$. Since $\eta(\chi\psi)<p$, it follows
that $\eta(\gamma^G)=1$, i.e $\gamma^G$ is a multiple of an irreducible.
\end{proof}
\begin{proof}[Proof of Theorem C]
{\bf (i)} Observe that $[\chi^2,\chi]=[\chi,\chi\overline{\chi}]$. Observe also
that $\Ker(\chi \overline{\chi})={\bf Z}(\chi)$. Thus if $[\chi,\chi\overline{\chi}]\neq 0$, then $\Ker(\chi)\geq {\bf Z}(\chi)$. Therefore $\Ker(\chi)={\bf Z}(\chi)$. Since 
${\bf Z}(G/\Ker(\chi))={\bf Z}(\chi)/\Ker(\chi)$ and $G$ is a $p$-group, it follows
that $G=\Ker(\chi)$ and so $\chi=1_G$.

{\bf (ii)} Assume that $[\chi^2,\lambda]\neq 0$ for some $\lambda\in \Lin(G)$.
Since $p\neq 2$ and $G$ is a $p$-group, 
there exists some character $\beta\in \Lin(G)$ such that $\beta^2=\lambda$. 
Thus $[\chi^2,\lambda]=[\chi^2, \beta^2]=[\chi\overline{\beta}, \overline{\chi}\beta]$.
Since $\chi\overline{\beta}, \overline{\chi}\beta \in \Irr(G)$, it follows that 
$$\chi\overline{\beta}= \overline{\chi}\beta=\overline{\chi\overline{\beta}}.$$
Since $\chi\overline{\beta}\in \Irr(G)$ is a real character 
and $G$ is a $p$-group with 
$p\neq 2$, it follows that  $\chi\overline{\beta}=1_G$. Thus $\chi=\beta$ and 
$\chi(1)=1$. 
 
{\bf (iii)} Observe that $\Ker(\chi^2)\geq \Ker(\chi)$. Working with the group 
$G/\Ker(\chi)$, without lost of generality we may assume that $\Ker(\chi)=1$.
We may also assume that $\chi(1)>1$. We are going to use induction on the order
of $|G|$.
 
 Let $Z={\bf Z}(\chi)$ be the center of the character $\chi$. Let $Y/Z$ be a chief 
factor of $G$. Let $\zeta\in \Lin(Z)$ be the unique character of $Z$ such that
 $\chi_Z=\chi(1)\zeta$. Since $G$ is a $p$-group and $Y/Z$ is a chief
 factor, it follows that $\zeta$ extends to $Y$.  Let $\iota\in \Lin(Y)$ be an 
 extension of $\zeta$. Since $Z={\bf Z}(\chi)$ and $Y/Z$ is a chief factor of 
 a $p$-group, it follows that $\iota$ lies
 below $\chi$. Let $G_{\iota}$ be the stabilizer of $\iota$ in 
 $G$. Observe that the stabilizer $G_{\iota^2}$ of $\iota^2$ contains 
 $G_{\iota}$. Since $|G:G_{\iota}|=p$, either $G=G_{\iota^2}$ or
 $G_{\iota^2}=G_{\iota}$. If $\eta(\chi^2)<p$, by Theorem A (i) 
 we have that $G_{\iota^2} = G_{\iota}$.
 If $G$ is a $p$-group with $p\neq 2$, then  $(\iota^g(y))^2=(\iota(y))^2$
 implies that $\iota^g(y)=\iota(y)$. Thus $G_{\iota^2}=G_{\iota}$ if
 $\eta(\chi^2)<p$ or $p\neq 2$.  

Let $\chi_{\iota}\in \Irr(G_{\iota})$ be the Clifford correspondent of 
$\chi$ and $\iota$, i.e. $\chi_{\iota}^G=\chi$ and $\chi_{\iota}$ lies above
$\iota$. Observe that $(\chi_{\iota}^2)_Y= \chi_{\iota}^2(1) \iota^2$ and 
$\iota^2\in \Lin(Y)$.
Thus all the irreducible constituents of the character 
$\chi_{\iota}^2$ lie above $\iota^2$. Since $G_{\iota}=G_{\iota^2}$,
by Clifford theory all the irreducible constituents of $\chi_{\iota}^2$ 
induce irreducibly to $G$. Since $|G_{\iota}|<|G|$, by induction 
there exist a 
subgroup $H$ of $G_{\iota}$ and a character $\alpha \in \Lin (H)$ such that 
$\alpha^{G_{\iota}}=\chi_{\iota}$ and $(\alpha^2)^{G_{\iota}}\in \Irr(G_{\iota})$.
Since all the irreducible constituents of $\chi_{\iota}^G$ induce 
irreducibly, it follows that $(\alpha^2)^G\in \Irr(G)$. Since $\chi_{\iota}^G=\chi$,
it follows that $(\alpha)^G=\chi$. Since $[ (\chi^2)_H, \alpha^2]\neq 0$ and 
$(\alpha^2)^G\in \Irr(G)$, it follows that $(\alpha^2)^G$ is an irreducible
constituent of $\chi^2$. Therefore $\chi^2$ has an irreducible constituent of
degree $\chi(1)$. 
\end{proof}

\end{section}

{\bf Acknowledgment.} I would like to thank Professor Everett C. Dade 
for helpful discussions.

\end{document}